\documentclass[reqno]{amsart}     
\usepackage{enumerate,amsmath,amssymb}     
\usepackage{pb-diagram}  
\usepackage{pstricks,pst-node,pst-coil,pst-plot}   
     
\theoremstyle{plain}      
 
\newtheorem{thm}{Theorem}[section]     
\newtheorem{notation}[thm]{Notation}

\newtheorem{prop}[thm]{Proposition}     
     
\theoremstyle{remark}      
     
\theoremstyle{definition}

\def\al{{\alpha}}

\def\la{{\lambda}}

\def\ep{{\varepsilon}}

\def\phi{{\varphi}}

\DeclareMathAlphabet{\doba}{U}{msb}{m}{n}

\gdef\mS{\doba{S}}

\def\eref#1{{\rm (\ref{#1})}}   

\let\<\langle 
\let\>\rangle 

\def\Vzero{V\setminus\{0\}}

\let\ss\scriptstyle

\begin{document}
\title{Second eigenvalue of the Yamabe operator and applications}

\maketitle     
\begin{center}
\sc S. El Sayed\footnote{elsayed@iecn.u-nancy.fr}
\end{center}
\begin{abstract}
Let $(M,g)$ be a compact Riemannian manifold of dimension $n\geq 3$. 
In this paper, we give various properties of the eigenvalues of the Yamabe operator $L_g$. In particular, we show  
how the second eigenvalue of $L_g$ 
is related to the existence of nodal solutions of the equation $L_g u = \ep \lvert u\lvert
^{N-2}u,$ where $\ep = +1,$ $0,$ or $-1.$
\end{abstract}
\begin{center}
\today
\end{center}

\tableofcontents

\section{Introduction}
This paper is part of a Phd thesis whose purpose is to study the relationships between the eigenvalues of the Yamabe operator,
 in particular their sign, and analytic, geometrical or topological properties of compact manifolds 
of dimension $n \geq 3$:\\

\noindent Let $(M,g)$ be a $n$-dimensional compact Riemannian manifold ($n \geq 3$). The {\em Yamabe operator} or {\em conformal Laplacian operator} $L_g$ is defined by
$$L_g(u)
:=c_n\Delta_gu+S_g u,$$
where $\Delta_g$ is the Laplace-Beltrami operator, $c_n=\frac{4(n-1)}{n-2}$ and
$S_g$ the scalar
curvature of $g$.
The Yamabe operator $L_g$ has discrete spectrum
$$\rm spec(L_g)=\left\lbrace \lambda_1(g), \lambda_2(g),\cdots \right\rbrace,$$
where the eigenvalues are such that
$$\la_1(g) < \la_2(g) \leq \la_3(g) \leq \cdots \leq \la_k(g) \cdots \to +\infty.$$
The $i-$th eigenvalue $\la_i(g)$ is characterized  by
\begin{eqnarray} \label{lambdai} 
\lambda_i (g) = \inf_{V\in Gr_{i}({H_1^2(M)})} \sup_{v \in \Vzero}
\frac{\int_M v L_{g} v \ dv_g }{\int_M v^2   \ dv_g},
\end{eqnarray}
where $Gr_{i}({H_1^2(M)})$ stands for the set of all $i$-dimensional subspaces of
$H_1^2(M).$ \\

\noindent Our project is to understand what we can deduce 
from the sign of $\lambda_i$. Now, we summarize what is known about this question and explain our motivations. At first, 
it is straightforward to see that the sign of $\lambda_1(g)$ is the same as the sign of the Yamabe constant $\mu(M,g)$ of 
$(M,g)$ (and as a consequence is conformally invariant). See Section \ref{Yamabesection} for more informations. 
Hence the positivity of $\lambda_1(g)$ has many consequences usually stated in terms of positivity of the Yamabe constant. For instance, we obtain

\begin{prop} \label{l1scal}
 A compact manifold $M$ of dimension $n \geq 3$ carries a metric with positive scalar curvature if and only if it carries a metric $g$ such that $\lambda_1(g) >0$.
\end{prop}
 We recall that classifying such compact manifolds is a challenging open problem, only solved for $n=3$ using Perelman's techniques. We also mention \cite{BarDahl} where  M. Dahl and C. B\"ar deduce many topological properties of compact manifolds from a careful study of the eigenvalues $\lambda_i$ of the Yamabe operator $L_g$.  \\ 

\noindent The sign of $\lambda_1$ can also be read in terms of existence or non-existence of positive solutions of the Yamabe equation: 
\begin{eqnarray} \label{yamabe_equation}
 L_g u = \epsilon |u|^{N-2} u,
\end{eqnarray}
where $N := \frac{2n}{n-2}$ and $\epsilon \in \{-1;0;1\}$. Inspired by this observation, B. Ammann and E. Humbert \cite{AmmannHumbert2005} enlighted the role of $\lambda_2$ in the existence of nodal solutions (i.e. having a changing sign) of the Yamabe equation \eref{yamabe_equation}. See again Section \ref{Yamabesection} for more explanations. \\

\noindent In this paper, we establish various properties of the eigenvalues of the Yamabe operator. First of all, we extend their definition to what we call  generalized metrics when possible (see Paragraph \ref{definitionsection}) and prove that their sign is a conformal invariant (see Paragraph \ref{conformalsection}). This paper initiates the study of the relationships between these conformal invariants and the topology of the manifold by showing that their negativity is not topologically obstructed (see Paragraph \ref{negativesection}). These investigations will be treated much more deeply in \cite{Elsayed2012}. The main point of this article is to complete the results of B. Ammann and E. Humbert \cite{AmmannHumbert2005} and to study how the sign of the second eigenvalue of the Yamabe operator can be related to the existence of nodal solutions of the Yamabe equation \eref{yamabe_equation}, in particular when the Yamabe constant of $(M,g)$ is negative. Our main result is to prove that under this condition, such a solution always exists with $\epsilon = sign (\lambda_2(g))$. This is the object of Theorem \ref{th1}.\\\\
The author would like to thank Emmanuel Humbert for his support and encouragements.

\section{Eigenvalues in conformal metrics}\label{definitionsection} 
In the whole paper, we will deal with the behavior of the eigenvalues of the Yamabe operator in a fixed conformal class. It will be usefull to express their definition relatively to a fixed metric. This is the goal of this section.
\subsection{Smooth metrics}
Let $(M,g)$ be a compact Riemannian manifold of dimension $n\geq 3$, we keep the notations of the introduction and for any metric $\widetilde g$, we will denote by 
$$\la_1(\widetilde g) < \la_2(\widetilde g) \leq \la_3(\widetilde g) \leq \cdots \leq \la_k(\widetilde g) \cdots \to +\infty,$$
 the eigenvalues of the Yamabe operator. We will deal with the case where $\widetilde g$ is conformal to $g$, i.e. when 
 $\widetilde g =
u^{N-2} g,$ where $u$ is a positive function of class $C^\infty$.
By referring to \cite{AmmannHumbert2005}, one sees that the $i-$th eigenvalue $\la_i(\widetilde
g)$ is given by
\begin{equation}\label{lam}
\lambda_i (\widetilde g) = \inf_{V\in Gr_{i}({H_1^2(M)})} \sup_{v \in \Vzero}
\frac{\int_M c_n |v|^2 + S_g v^2 \, dv_g }{\int_M v^2 \ u^{N-2}  \ dv_g},
\end{equation}
where $\widetilde g = u^{N-2}g;$ $u\in C^\infty(M),$ $u>0$ and $Gr_{i}({H_1^2(M)})$ stands for the set of all $i$-dimensional subspaces of
$H_1^2(M).$
\subsection{Generalized metrics}\label{negativesection}
 Reducing to smooth metrics will too restrictive for our investigations. We will need to work with {\em generalized metrics}, i.e. metrics of the form  $\widetilde g = u^{N-2}g$ with $u\in L^N(M),$ $u\geq 0$ and $u\not\equiv 0$. 
The Yamabe operator $L_{\widetilde g}$ has no meaning any more but the definition of $\lambda_i(\widetilde g)$ can anyway be extended to this case by using \eref{lam} as it was done in \cite{AmmannHumbert2005} when the Yamabe constant was non-negative i.e. when $\lambda_1(g) \geq 0$.
When $\lambda_1(g) < 0$, the situation is a little bit different: $\lambda_i(\widetilde g)$ defined by \eref{lam} can be $-\infty$ as proved by the following proposition.
\begin{prop} \label{lainfty}
 Assume that $\lambda_1(g)<0$, then there exists $u\in L^N(M)$, $u \not\equiv0$, $u\geq 0$ such that $\lambda_1(\widetilde g) = -\infty$,  where $\widetilde g = u^{N-2}g$.
\end{prop}
\noindent This proposition will be proved in Paragraph \ref{infini}. To make sure that $\lambda_1(\widetilde g)$ is finite, one has to assume in addition that $u$ is positive.
\begin{prop}\label{la_1bornee}
Let $u$ be a positive function in $L^N(M)$. Suppose that $\la_1(g)<0$. Then, we have 
$$\la_1(\widetilde g)> -\infty.$$
\end{prop}
The proposition is proved in Paragraph \ref{preuveprop}.
\begin{notation}
 The $i^{th}$ eigenvalue of $L_g,$ $\la_i(\widetilde g) = \la_i(u^{N-2}g)$ will be
denoted by $\la_i(u)$ when there is no ambiguity about $g$. 
\end{notation}
\subsubsection{Proof of Proposition \ref{lainfty}}\label{infini}
We have $\la_1(g)<0$, this implies that there exists a function $v\in C^{\infty}(M)$ such that 
$$\int_M (L_g v)v\,dv_g < 0.$$
Let $P$ be a point of $M$. For $\ep>0$, we define $\eta_{\ep}$ as follows
$$\left\{
\begin{array}{c}
0 \leq \eta_{\ep} \leq 1,\\\\
\eta_{\ep} = 0 \text{ on }B_{\ep}(P),\\\\
\eta_{\ep} = 1\text{ on } M\backslash B_{2\ep}(P),\\\\
\lvert\nabla \eta_{\ep} \lvert\leq \frac{2}{\ep}.
\end{array}\right.$$
where $B_{\delta}(P)$ stands for the ball of center $P$ and radius $\delta$ in the metric $g$.
Then one easily cheks
$$\lim_{\ep\rightarrow 0}\int_M (L_g(\eta_{\ep}v))(\eta_{\ep}v)\,dv_g = \int_M (L_g v)v\,dv_g.$$
We define $w:= \eta_{\ep}v$. Therefore, for a fixed small $\ep > 0$, we have
$$\int_M (L_g w)w\, dv_g < 0.$$
Let $u\geq 0$, $u\not\equiv 0$ of class $C^{\infty}$ with support in $B_{\ep}(P).$ For $\alpha >0$, since $(w+\alpha)u \not\equiv 0$,
we can write 
\begin{eqnarray*}
 \la_1(\widetilde g) &=& \inf_{v'} \frac{\int_M (L_g v')v'\,dv_g}{\int_M u^{N-2} v'^2 \,dv_g}\\
&\leq& \lim_{\alpha \rightarrow 0^+} \frac{\int_M (L_g (w+\alpha))(w+\alpha)\,dv_g}{\int_M u^{N-2} (w+\alpha)^2 \,dv_g}. 
\end{eqnarray*}
Moreover, we have 
$$\lim_{\alpha\to 0^+}\int_M (L_g (w+\alpha))(w+\alpha)\,dv_g=\int_M (L_g (w))w\,dv_g <0,$$
and $$ \lim_{\alpha\to 0^+} \int_M u^{N-2} (w+\alpha)^2 \,dv_g = 0$$
which gives that 
$$\lim_{\alpha\rightarrow 0}\frac{\int_M (L_g (w+\alpha))(w+\alpha)\,dv_g}{\int_M u^{N-2} (w+\alpha)^2 \,dv_g} = -\infty.$$
This ends the proof of Proposition \ref{lainfty}.
\subsubsection{Proof of Proposition \ref{la_1bornee}}\label{preuveprop}
 Let $(v_m)_m$ be a minimizing sequence for $\la_1(u),$ i.e. $v_m\in
H_1^2(M)$ such that
$$\lim_{m\longrightarrow \infty} \frac{\int_M c_n\lvert \nabla v_m\lvert^2 + S_g
v_m^2 \ dv_g}{\int_M \lvert u\lvert^{N-2} v_m^2 \ dv_g} = \lim_{m\longrightarrow \infty} \la_m = \la_1(u) < 0.$$
Since $(\lvert v_m\lvert)_m$ is also a minimizing sequence for $\la_1(u),$ we can
assume that $v_m\geq 0.$ We normalize $v_m$ by $\int_M \lvert u\lvert^{N-2} v_m^2 \
dv_g = 1.$ Here we show that $(v_m)_m$ is bounded in $H_1^2(M).$ Indeed, suppose that $(v_m)_m$
is not bounded in $H_1^2(M)$ and let 
$$v_m^\prime = \frac{v_m}{\parallel v_m\parallel_{H_1^2(M)}}.$$ 
$(v_m^\prime)_m$ is bounded in $H_1^2(M),$ and his norm is equal to 1, then there
exists $v^\prime \in H_1^2(M),$ (after restriction to a subsequence) such that 
$$v_m^\prime \rightharpoonup v^\prime \text{ in } H_1^2(M),$$
$$v_m^\prime\longrightarrow v^\prime \text{ in }L^2(M).$$
We have 
$$c_n \int_M \left|\nabla v_m'\right|^2\, dv_g + \int_M S_g v_m'^2 \, dv_g = \la_m \int_M \left|u\right|^{N-2}v_m'^2\, dv_g.$$
Moreover,
$$\int_M \left|u\right|^{N-2}v'^2\, dv_g \leq \int_M \left|u\right|^{N-2}v_m'^2\, dv_g \to_{m\longrightarrow \infty} 0$$
since 
$$\left\|v_m\right\|_{H_1^2(M)}\longrightarrow \infty.$$
It follows that  $$\int_M \left|u\right|^{N-2}v'^2\, dv_g = 0$$
and since $u$ is positive,  
$$v' = 0.$$
Now, we write
$$1 = \int_M \left|\nabla v_m'\right|^2\, dv_g + \underbrace{\int_M \left|v_m'\right|^2 \, dv_g}_{\longrightarrow 0}.$$
We deduce that $$\lim_{m\to \infty} \int_M \left|\nabla v_m'\right|^2\, dv_g=1,$$
giving the desired contradiction: 
$$c_n \underbrace{\int_M \left|\nabla v_m'\right|^2\, dv_g}_{\longrightarrow 1} + \underbrace{\int_M S_g v_m'^2 \, dv_g}_{\longrightarrow 0} = \la_m \int_M \left|u\right|^{N-2}v_m'^2\, dv_g \leq 0.$$
This proves that $(v_m)_m$ is bounded in $H_1^2(M),$ and implies that $\la_m \geq C.$ We finally get $\la_1(\widetilde g) > -\infty.$

\subsection{PDE associated to $\la_i$}\label{}
\begin{prop}\label{existancedevetw}
For any non-negative function $u\in L^N(M),$ such that $\la_1(u)>-\infty$, there exists 
functions $v_1>0, v_2, \ldots , v_k \in H_1^2(M)$ having a changing sign, such that in the sense of
distributions, we have
$$L_g v_1 = \la_1(u) \lvert u\lvert^{N-2} v_1,$$
and $$L_g v_k = \la_k(u) \lvert u\lvert^{N-2} v_k.$$ 
Moreover, we can normalize the $v_k$ by $$ \int_M \lvert u\lvert^{N-2} v_k^2 \ dv_g = 1 \text{ and }\int_M
\lvert u\lvert^{N-2} v_i v_j \ dv_g =0 \hskip0.2cm \forall i \neq j.$$
\end{prop}
\textbf{Proof:} Let $(v_m)_m$ be a minimizing sequence for $\la_1(u),$ i.e. $v_m\in
H_1^2(M)$ such that
$$\lim_{m\longrightarrow \infty} \frac{\int_M c_n\lvert \nabla v_m\lvert^2 + S_g
v_m^2 \ dv_g}{\int_M \lvert u\lvert^{N-2} v_m^2 \ dv_g} = \la_1(u).$$
According to the Paragraph \ref{preuveprop}, we get that $(v_m)_m$ is bounded in $H_1^2(M)$ and there exists
$v\geq 0$ in $H_1^2(M)$ such that $v_m$ converges to $v$ weakly in $H_1^2(M)$ and
strongly in $L^2(M)$ (after restriction to a subsequence). We now want to prove 
\begin{eqnarray}\label{equa}
 \int_M \lvert u\lvert^{N-2} v^2 \ dv_g = \lim_{m\longrightarrow \infty} \int_M
\lvert u\lvert^{N-2} v_m^2 \ dv_g = 1.
\end{eqnarray}
If $u$ is smooth, this relation is clear. So let us assume that $u\in \L^N(M),$ let $A$ be a large real number and set $u_A = \inf\left\lbrace
u,A\right\rbrace.$ By H\"{o}lder inequality, we write
\begin{eqnarray*}
 \left|  \int_M u^{N-2} \left(v_m^2 - v^2\right) \,dv_g \right| &=&  \left|  \int_M
\left(u^{N-2} -u_A^{N-2} +u_A^{N-2}\right)\left(v_m^2 - v^2\right) \,dv_g \right|\\
& \leq &  \left( \int_M u_A^{N-2} |v_m^2 - v^2| \,dv_g + \int_M (u^{N-2}
       - u_A^{N-2}) (|v_m|+|v|)^2 \,dv_g\right) \\
        & \leq & A^{N-2} \int_M |v_m^2 - v^2| \,dv_g\\ 
        &&+ {\left( \int_M (u^{N-2}-u_A^{N-2})^\frac{N}{N-2} \,dv_g
\right)}^{\frac{N-2}{N}}  {\left( \int_M (|v_m|+|v|)^N \,dv_g
  \right)}^{\frac{2}{N}}.
\end{eqnarray*}
$(v_m)_m$ is bounded in $H_1^2(M),$ it is bounded in $L^N(M).$ Hence there exists a
constant $C$ such that 
$$\int_M \left(|v_m|+|v|\right)^N \,dv_g \leq C.$$
The convergence in $L^2(M)$ gives
$$\lim_{m\longrightarrow \infty}\int_M |v_m^2 - v^2| \ dv_g = 0.$$
By dominated convergence theorem, we have
$$\lim_{A\longrightarrow \infty}\int_M \left(u^{N-2}-u_A^{N-2}\right)^\frac{N}{N-2}
\ dv_g= 0.$$
Hence, we get 
\eref{equa}.
Since $$\lim_m \int_M \langle\nabla v_m, \nabla \phi\rangle \ dv_g = \int_M
\langle\nabla v, \nabla \phi\rangle \ dv_g,$$
$$\lim_m \int_M S_g v_m \phi \ dv_g = \int_M S_g v \phi \ dv_g$$
and 
$$\lim_m \int_M \lvert u\lvert^{N-2} v_m \phi \ dv_g = \int_M \lvert u\lvert^{N-2} v
\phi \ dv_g,$$
(by strong convergence in $L^2(M)$),
we obtain that in the sense of distributions $v$ verifies 
$$L_gv = \la_1(u) \lvert u\lvert^{N-2} v.$$
Now we define
$$\la_k^\prime (u) = \inf_{{\ss v_k; \lvert u\lvert^{\frac{N-2}{2}} v_k\not\equiv 0\atop
\ss \int_M \lvert u\lvert^{N-2} v_i v_k \ dv_g = 0 \forall i<k}} \frac{\int_M c_n\lvert \nabla
v_k\lvert^2 + S_g v_k^2 \ dv_g}{\int_M \lvert u\lvert^{N-2} \lvert v_k\lvert^2 \ dv_g}.$$
we remark that $\la_k^\prime (u) = \la_k(u)$ and $v_k$ is constructed by induction using the same method. This ends the proof of Proposition \ref{existancedevetw}.
\hfill $\square$
\section{Sign of $\la_i$}
\subsection{The sign of $\la_i$ is conformally invariant} \label{conformalsection}
\begin{prop}\label{si}
The sign of $\la_i$ is independent of the metric selected in the conformal class. More precisely, for any conformal metric $\widetilde g = u^{N-2}g,$ where $u$ is a non-negative function in $L^N(M)$, $\la_i(u)$ and $\la_i(1)$ have same sign.
\end{prop}
\textbf{Proof:} We assume for example that $\la_i(u) = 0$ and $\la_i(1) > 0,$ we
know that
 $$\la_i (u) = \inf _{u_1,\ldots, u_i}\sup_{\la_1,\ldots,\la_i}\frac{\int_M
L_g(\la_1 u_1+ \cdots+\la_i u_i)(\la_1 u_1+ \cdots+\la_i u_i) \ dv_g}{\int_M (\la_1
u_1+ \cdots+\la_i u_i)^2 u^{N-2} \ dv_g},$$
and $$\la_i (1) = \inf _{u_1,\ldots, u_i}\sup_{\la_1,\ldots,\la_i}\frac{\int_M
L_g(\la_1 u_1+ \cdots+\la_i u_i)(\la_1 u_1+ \cdots+\la_i u_i) \ dv_g}{\int_M (\la_1
u_1+ \cdots+\la_i u_i)^2 \ dv_g}.$$
Suppose that $\la_i(u)$ is attained by $v_1,\ldots,v_i$. Since denominators of this expressions are positive, then
$$\sup_{\la_1,\cdots,\la_i}\int_M L_g(\la_1 v_1+ \cdots+\la_i v_i)(\la_1 v_1+
\cdots+\la_i v_i) \ dv_g = 0.$$
So
$$\la_i(1) \leq \sup_{\la_1,\ldots,\la_i}\frac{\int_M L_g(\la_1 v_1+ \cdots+\la_i
v_i)(\la_1 v_1+ \cdots+\la_i v_i) \ dv_g}{\int_M (\la_1 v_1+ \cdots+\la_i v_i)^2
u^{N-2} \ dv_g} = 0,$$
which gives a contradiction.
The remaining cases are treated similarly.
\subsection{The negativity of $\la_k$ is not topologically obstructed}
In this paragraph, we will see that on each manifold, there exists a metric which has a negative $k^{th}$-eigenvalue.
\begin{prop}\label{lak}
On any compact Riemannian manifold $M,$ and for all $k\geq 1$ there exists a metric
$g$ such that
$$\la_k(g)<0.$$
\end{prop}
\textbf{Proof:} Let $M$ be a compact Riemannian manifold of dimension $n,$ and we
take $k$ spheres of dimension $n=dim (M).$ We equip each sphere $\mathbb{S}^n$ by the
same metric $g,$ such that $\mu(g)<0.$ We can do this by referring to \cite{Aubin98}
(Theorem [1] page 38).
Let $P\in M,$ since $\mu(g)<0,$ for all $\ep, \delta >0$ we can find a function $u$
supported in
$\mathbb{S}^n\backslash
{B_{\ep}(P)}$ such that
$$\frac{\int_M c_n\lvert\nabla u\lvert^2 + S_g u^2 \ dv_g}{\int_M u^2 \ dv_g}<
-\delta.$$
Indeed, let $\eta_\ep$ be a smooth cut-off function such that $0\leq\eta_\ep\leq 1,$
$\eta_\ep
(B_{\ep}(P))=0,$ $\eta_\ep(\mathbb{S}^n\backslash
{B_{2\ep}(P)})=1,$ $\lvert\nabla \eta_\ep\lvert\leq \frac{2}{\ep}$ and a function $v$
satisfying $$I_g(v) = \frac{\int_M
c_n\lvert\nabla v\lvert^2 + S_g v^2 \ dv_g}{\int_M v^2 \ dv_g}< -2\delta.$$
Note that the existence of $v$ is given by the fact that $\mu(g)<0.$ The desired
function $u$
will be given by $\eta_\ep v,$ where $\ep>0$ is sufficiently small. Indeed, it
suffices to notice that, as easily checked,
$$\lim_{\ep\longrightarrow 0}I_g(\eta_\ep v) = I_g(v).$$
Let $P_1,\cdots,P_k$ be points of $M.$ We consider the following connected sum
$$M^\prime = M \# (\mathbb{S}^n)_1 \# \ldots \# (\mathbb{S}^n)_k,$$
where the $(\mathbb{S}^n)_i$ are attached at $P$ on the spheres $\mathbb{S}^n$ and at
$P_i$ on $M$ so that the handles are attached in $B_{\ep}(P)$ and $B_{\ep}(P_i)$. Note
that
$M^\prime$ is diffeomorphic to $M.$ Moreover, the above construction allows to
see $(\mathbb{S}^n)_i\backslash
B_{\ep}(P_i)$ as a part of $M^\prime.$\\
We take on $M^\prime$ any metric $h$ satisfying
$$h\arrowvert_{(\mathbb{S}^n)_i\backslash B_{\ep}(P_i)}= g.$$
On $M^\prime,$ we define the following function
\begin{eqnarray*}
   u_i=\left|\;
   \begin{matrix}
     u\hfill & \hbox{on } (\mathbb{S}^n)_i\backslash B_{\ep}(P_i)\\\\
     0\hfill
     & \hbox{otherwise} .
   \end{matrix}
   \right.
\end{eqnarray*}
Since the $u_i$ have disjoint supports, we get
\begin{eqnarray*}
 \la_k(1)&\leq& \sup_{\la_1,\ldots,\la_k}\frac{\int_{M^\prime} L_g(\la_1 u_1+
\cdots+\la_k
u_k)(\la_1 u_1+ \cdots+\la_k u_k) \ dv_h}{\int_{M^\prime} (\la_1 u_1+ \cdots+\la_k
u_k)^2 \
dv_h}\\\\
&\leq&  \sup_{\la_1,\ldots,\la_k}\frac{(\la_1^2+ \ldots + \la_k^2)\int_M (L_gu)u \
dv_g}{(\la_1^2+ \ldots + \la_k^2)\int_M u^2 \ dv_g}\\\\
&\leq& \frac{\int_M (L_gu)u \ dv_g}{\int_M u^2 \ dv_g}< -\delta.
\end{eqnarray*}

\section{Nodal solutions of the Yamabe equations}\label{Yamabesection}
\noindent A famous problem in Riemannian geometry is the Yamabe problem, solved between 1960 and 1984 by Yamabe, Tr\"udinger, Aubin and Schoen, \cite{yamabe:60, trudinger:68, Aubin1976, Schoen:84}. The reader can also refer to \cite{lee.parker:87, Hebey97, Aubin98}. The Yamabe problem consists in finding a metric $\widetilde g$ conformal to 
$g$ such that the scalar curvature $S_{\widetilde g}$ of $\widetilde g$ is constant.
Solving this problem is equivalent to finding a positive smooth function and a
number $C_0 \in \mathbb{R}$ such that \begin{equation}\label{muun}
L_g(u) = C_0 \lvert u\lvert^{N-2} u,\end{equation}
where $N = \frac{2n}{n-2}.$ In order to obtain solutions of the Yamabe equation we
define the Yamabe invariant by $$\mu(M,g):= \inf_ {u\neq 0, u\in C^\infty
(M)}Y(u),$$
where $$Y(u) = \frac{\int_M c_n\lvert \nabla u\lvert^2+S_g u^2 \ dv_g}{\left(\int_M \lvert
u\lvert^N \ dv_g\right)^\frac{2}{N}}.$$ 
The works of Yamabe, Tr\"udinger, Aubin and Schoen provides a positive smooth minimizer $u$ of $Y$, satisfying, if normalized by
$\|u\|_{L^N(M)} = 1,$
$$L_gu = \mu(M,g) \lvert u\lvert^{N-2}u.$$
The metric $\widetilde g = u^{N-2} g$ is the desired metric: its scalar curvature is constant equal to $\mu(M,g)$.
If we set $u' = \mu(M,g)^{\frac{n-2}{4}} u$, we obtain a positive solution of 
$$L_gu' = \ep \lvert u'\lvert^{N-2}u'$$
where $\ep = \hbox{ sign } (\mu(M,g))= \hbox{ sign}(\lambda_1(g))$. \\

\noindent Now, if $\mu(M,g)\geq 0,$ it is easy to chek that 
$$\mu(M,g) = \inf_{\widetilde g\in \left[ g\right]} \la_1(\widetilde g) vol
(M,\widetilde g)^\frac{2}{n},$$
where $\left[ g\right]$ is the conformal class of $g$ and $\la_1$ is the first
eigenvalue of the Yamabe operator $L_g$.
Inspired by this approach,  in their paper
\cite{AmmannHumbert2005}, B. Ammann et E. Humbert introduced the second Yamabe invariant defined by
\begin{eqnarray*}
 \mu_2(M,g)&=& \inf_{\widetilde g} \la_2(\widetilde g) vol(M,\widetilde
g)^\frac{2}{n}\\
&=& \inf_u \la_2(u^{N-2}g) \left(\int_M u^N \ dv_g\right)^\frac{2}{n}.
\end{eqnarray*}
They studied this invariant in the case where $\mu(M,g)\geq 0$, and they proved that
$\mu_2$ is attained by a generalized metric, (i.e. a metric of the form $u^{N-2}g$
where $u\in L^N(M)$, $u\geq 0$ which may vanish), in the following two cases\\
$\bullet$ $\mu(M,g)>0$, $(M,g)$ is not locally conformally flat and $n\geq 11$.\\
$\bullet$ $\mu(M,g)=0$, $(M,g)$ is not locally conformally flat and $n\geq 9$.\\ 
In this context, they proved that $u$ is the absolute value of a changing sign function $w$ of
class $C^{3,\alpha}(M)$, which verifies the following
equation
$$L_g w = \mu_2(M,g)\lvert w\lvert^{N-2} w.$$
Many works are devoted to the study of this kind of solutions, for example \cite{AmmannHumbert2005}, 
\cite{djadli.jourdain:02}, \cite{holcman:99}, \cite{hebeyvaugon:94}, \cite{Vetois2007}. See also \cite{BenaliliBoughazi2008} 
for an analogue study for the Paneitz-Branson operator. 
Setting again $w' = \mu_2(M,g)^{\frac{n-2}{4}} w$, we obtain a solution of 
$$L_g w'= \ep\lvert w'\lvert^{N-2} w'$$
with $\epsilon=1= \hbox{ sign }(\mu_2(M,g))= \hbox{ sign }(\lambda_2(g)).$
The goal of this section is to study if this result extends to metrics where the sign of $\la_2(g)$ is arbitrary.\\

\noindent The answer is yes when $\la_2 < 0$ without any other condition, we obtain this result by a method different than the one of \cite{AmmannHumbert2005}. Notice that this situation occurs for a large number of metrics (see Proposition \ref{lak}).                                                                  
When $\la_2\geq 0$, we show that the methods in \cite{AmmannHumbert2005} can be extended to the case where 
$\mu(M,g)<0$.  Namely, the main result of this paper is:  
\begin{thm}\label{th1}
Let $(M,g)$ be a compact Riemannian manifold of dimension $n\geq3$ whose Yamabe
invariant $\mu(M,g)$ is strictly negative, we denote by $\la_2$ the second
eigenvalue of
$L_g.$ Then, if $\la_2\leq 0$ or if $\la_2>0$, $(M,g)$ not locally conformally flat and $n\geq 6$:\\
There exists a function $w$ changing sign, solution of the equation
$$L_g w = \ep \lvert w \lvert^{N-2}w,$$
where $\ep=+1$ if $\la_2>0,$ $\ep = -1$ if $\la_2< 0$ and $\ep=0$ if $\la_2 = 0.$
Moreover, $w \in C^{3,\alpha}(M)$, for all $\alpha< N-2.$
\end{thm}
\subsection{The case $\la_2 = 0$}
This case is obvious: indeed, Proposition \ref{existancedevetw} provides the existence of a nodal solution 
$v$ of  
$L_g v = 0= \ep  \left|v\right|^{N-2}v$ 
where $\ep = 0 = \hbox{ sign }(\lambda_2(g))$. 
\subsection{The case $\la_2 > 0$}\label{section}
As in \cite{AmmannHumbert2005}, we introduce
 the second Yamabe invariant given
by
\begin{eqnarray*}
 \mu_2(M,g)&=& \inf_{\widetilde g}\la_2(\widetilde g) vol(M,\widetilde g)^\frac{2}{n}\\
&=& \inf_{u>0} \la_2(u)\left(\int_M u^N \ dv_g\right)^\frac{2}{n}.
\end{eqnarray*}
By Proposition \ref{uw} below, the problem reduces to finding a minimizer of 
$\mu_2(M,g)$. The case where $\mu(M,g) \geq 0$ have been treated in \cite{AmmannHumbert2005}. We will then focus on the case where $\mu(M,g) < 0$ (i.e. $\la_1(g) <0$). 
We will see that the method of Ammann and Humbert remains valid in this case and the
following three  propositions answer our questions.\\

\begin{prop} \label{uw}
Let $(M,g)$ be a compact Riemannian manifold of dimension $n\geq3$, such that $\la_2>0$. If
\begin{equation}\label{in}
\mu_2(M,g)<\mu(\mathbb{S}^n),
\end{equation}
with
$\mu(\mathbb{S}^n)=n(n-1)\omega_n^{\frac{2}{n}}$, where $\omega_n$ stands for the
volume of the standard sphere $\mathbb{S}^n$, then the second Yamabe invariant is attained by a non-negative function $u\in L^N(M)$
that we normalize by $\int_M u^N \ dv_g = 1$. There exists a
function $w$ having a changing sign which verifies in the sense of distributions the
following equation
\begin{equation}
 \label{mu2}
L_gw = \mu_2(M,g) \lvert u\lvert^{N-2} w.
\end{equation}
The functions $u$ and $w$ will be normalized by 
$$\int_M u^N \ dv_g = 1, \hskip0.2cm
\int_M u^{N-2} \ w^2 \ dv_g = 1.$$
\end{prop}
\begin{prop}\label{a}
The two functions $u$ and $w$ given by Proposition \ref{uw} satisfy 
$$u = \lvert w\lvert.$$
\end{prop}
Finally, we give a condition under which assumption (\ref{in}) is satisfied:
\begin{prop}\label{ip}
Let $(M,g)$ be a compact Riemannian manifold of dimension $n\geq 6,$ suppose that
$M$ is not locally conformally flat and his Yamabe invariant $\mu(M,g)<0,$ then
$$\mu_2(M,g)< \mu(\mathbb{S}^n).$$
\end{prop}
\subsubsection{Proof of Proposition \ref{uw}}
The case where $\mu(M,g)  \geq 0$ is done in \cite{AmmannHumbert2005}, hence we consider here the case where $\mu(M,g)  < 0$. By the solution of the Yamabe problem, we can assume without loss of generality, that 
$S_g = -1$.
Let $(u_m)_m$ be  a minimizing sequence for $\mu_2(M,g),$ i.e., $u_m$ is positive, smooth and 
$$\lim_{m\longrightarrow
\infty}\lambda_2(u_m)\left(\int_Mu_m^Ndv_g\right)^\frac{2}{n} =
\mu_2(M,g).$$
The sequence $(u_m)_m$ will be choosen such that $\int_M u_m^N \ dv_g=1$, hence
$\mu_2(M,g)=\lim_{m\longrightarrow \infty}\lambda_2(u_m).$ For each $u_m,$ Proposition \ref{existancedevetw} provides  the existence of a function $w_m \in H_1^2(M)$ such that
\begin{equation}\label{eq2}
L_g w_m = \la_2(u_m)\lvert u_m\lvert^{N-2}w_m.
\end{equation}
Moreover, the sequence $(w_m)_m$ can be normalized by 
$$\int_M\lvert u_m\lvert^{N-2}w_m^2 \ dv_g=1.$$
Since $\int_Mu_m^Ndv_g=1,$ $(u_m)_m$ is bounded in $L^N(M)$ which is a
reflexive space, there exists $u\in L^N(M)$ such that $u_m$ converges weakly to
$u$ in $L^N(M)$, we have
$$u_m\rightharpoonup u \text{ in }L^N(M).$$
$\bullet$ The sequence $(w_m)_m$ is bounded in $H_1^2(M).$\\
We proceed by contradiction and assume that 
$\|w_m\|_{H_1^2(M)}\longrightarrow \infty.$ Let 
$$w^\prime_m=\frac{w_m}{\|w_m\|_{H_1^2(M)}}.$$
$\|w^\prime _m\|_{H_1^2(M)}=1,$ hence $(w^\prime _m)_m$ is bounded in
$H_1^2(M).$ Since $H_1^2(M)$ is a reflexive space, this implies using Kondrakov and Banach-Alaoglu theorems, that there
exists a subsequence $(w_m^\prime)_m$ and $w^\prime \in H_1^2(M)$ such that
$$w^\prime_m\rightharpoonup w^\prime \text{ in } H_1^2(M),$$
and
$$w^\prime_m\longrightarrow w^\prime \text{ in } L^2(M).$$
Equation (\ref{eq2}) is linear, so $w^\prime _m$ satisfies
$$L_gw_m^{\prime} = \lambda_2(u_m)\left|u_m \right|^{N-2}
w^\prime _m.$$
Hence for all $\varphi \in C^\infty (M),$ we have:
$$c_n \int _M \left\langle \nabla w^\prime _m,\nabla\varphi\right\rangle
dv_g+\int_MS_gw^\prime _m\varphi dv_g=\int_M 
\lambda_2(u_m)\left|u_m\right|^{N-2}w^\prime_m \varphi dv_g.$$
Since $w^\prime _m \rightharpoonup w^\prime$ in $H_1^2(M)$ and $w \longmapsto
\left\langle \nabla w ,\nabla\varphi\right\rangle$ is a linear form on $H_1^2(M),$
then
$$c_n\int_M \left\langle \nabla w_m^\prime,\nabla\varphi  \right\rangle dv_g
\longrightarrow c_n \int_M\left\langle \nabla  w^{\prime},\nabla\varphi\right\rangle
dv_g.$$
The sequence $w_m^\prime$ converges strongly to $w^\prime$ in $L^2(M)$. This gives that $$\int_M S_gw_m^\prime \varphi
dv_g\longrightarrow \int_M S_g w^\prime \varphi dv_g.$$
Using H\"older inequality, we obtain that
$\int_M\left|u_m\right|^{N-2}w^\prime_m\varphi \ dv_g\longrightarrow 0.$ Indeed, 
\begin{eqnarray*}
\left|\int_M\left|u_m\right|^{N-2}w^\prime_m\varphi
dv_g\right|&\leq&\left\|\varphi\right\|_{\infty}\int_M\left|u_m\right|^\frac{N-2}{2}\left|w_m^\prime\right|\left|u_m\right|^\frac{N-2}{2}
\ dv_g\\
&\leq&
\left\|\varphi\right\|_{\infty}\left(\int_M\left|u_m\right|^{N-2}{w^\prime_m}^2 \
dv_g\right)^\frac{1}{2}\left(\int_M\left|u_m\right|^{N-2} \ dv_g\right)^\frac{1}{2}\\
&\leq&\left\|\varphi\right\|_{\infty}\frac{\left(\int_M\left|u_m\right|^{N-2}{w_m}^2 \
dv_g\right)^\frac{1}{2}}{\left\|w_m\right\|_{H_1^2(M)}}\left(\int_M\left|u_m\right|^{N}
\
dv_g\right)^\frac{N-2}{2N}\left(\rm vol(M,g)^{1-\frac{N-2}{N}}\right)^\frac{1}{2}\\
&\leq&\left\|\varphi\right\|_{\infty} \frac{1}{\left\|w_m\right\|_{H_1^2(M)}}\left(\rm
vol(M,g)\right)^\frac{1}{N} \longrightarrow_{m\longrightarrow +\infty}0.
\end{eqnarray*}
Then
$$c_n \int_M\left\langle \nabla  w^{\prime},\nabla\varphi\right\rangle \ dv_g+\int_M
S_g w^\prime \varphi \ dv_g=0,$$
which means that in the sense of distributions, we have
$$L_g w^\prime=0.$$
Since $\la_1(1)<0$ and $\la_2(1)$ is
positive, $0\notin Sp(L_g)$. It follows that $w^\prime
=0.$ Now, we also have
$$\int_M c_n\left|\nabla w^\prime _m\right|^2 \ dv_g+\int_M S_g{w^\prime_m}^2 \
dv_g=\lambda_2(u_m)\int_M\left|u_m\right|^{N-2}{w_m^\prime}^2 \ dv_g,$$
with
$$\lambda_2(u_m)\int_M\left|u_m\right|^{N-2}{w_m^\prime}^2 \
dv_g=\frac{\lambda_2(u_m)}{\left\|w_m\right\|^2_{H_1^2(M)}}\longrightarrow 0$$ and
$$\int_M S_g{w^\prime_m}^2 \ dv_g\longrightarrow \int_M S_g{w^\prime}^2 \ dv_g = 0.$$
Hence $$\int_M\left|\nabla w^\prime_m\right|^2 \ dv_g\longrightarrow 0.$$
Finally, we get that
$$\left\|w^\prime_m\right\|^2_{H_1^2(M)}=1=\int_M\left|\nabla w_m^\prime\right|^2 \
dv_g+\int_M{w^\prime_m}^2 \ dv_g\longrightarrow 0,$$
which gives the desired contradiction. We obtain that $(w_m)_m$ is a bounded
sequence in $H_1^2(M).$ Then there exists $w\in H_1^2(M)$ such that:
$$w_m\rightharpoonup w \text{ in } H_1^2(M),$$
$$w_m\longrightarrow w \text{ in }L^2(M).$$
It follows that in the sense of distributions, we have 
$$L_gw=\mu_2(M,g)\left|u\right|^{N-2}w.$$
It remains to show that $w$ changes sign and is different from zero.\\\\
$\bullet$ Suppose that $w$ does not change sign. Without loss of generality, we
can assume that $w\geq 0.$ In the sense of distributions, we have
\begin{equation}\label{eg}
c_n\Delta_gw+S_gw=\mu_2(M,g)\left|u\right|^{N-2}w.\end{equation}
It was already mentioned at the beginning of this section that we can assume that
$S_g<0,$ because $\mu(M,g)<0.$ Integrating (\ref{eg}) over $M$, we get:
$$\underbrace{\int_Mc_n\Delta_gw \ dv_g}_{=0}+\underbrace{\int_MS_gw \
dv_g}_{<0}=\underbrace{\mu_2(M,g)\int_M\left|u\right|^{N-2}w \ dv_g}_{\geq 0}.$$
This gives a contradiction unless $w \equiv 0$ which is prohibited by what follows. \\

$\bullet$ Assume that $w=0.$ By referring to \cite{Hebey97} and \cite{Aubin1976} we
have the following theorem:\\
If $(M,g)$ is a Riemannian manifold of dimension $n\geq
3,$ for all
$\epsilon >0,$ there exists $B_{\epsilon}$ such that for any $u\in H_1^2(M),$ we have
$$\left(\int_M\left|u\right|^N \ dv_g \right)^\frac{2}{N}\leq
(\mu(\mathbb{S}^n)^{-1}+\epsilon)\left(\int_M c_n\left|\nabla u\right|^2 \
dv_g+B_\epsilon\int_Mu^2 \ dv_g\right).$$
We obtain
\begin{eqnarray*}
c_n\int_M\left|\nabla w_m\right|^2 \ dv_g+S_g\int_Mw_m^2 \
dv_g&=&\mu_2(M,g)\int_M\left|u_m\right|^{N-2}w_m^2 \ dv_g\\
&\leq&\mu_2(M,g)\underbrace{\left(\int_M\left|u_m\right|^{N} \
dv_g\right)^\frac{N-2}{N}}_{=1}\left(\int_M\left|w_m\right|^{N} \
dv_g\right)^\frac{2}{N}\\
&\leq&\underbrace{\mu_2(M,g)(\mu(\mathbb{S}^n)^{-1}+\epsilon)}_{<1
(\text{if}\hskip0.1cm \ep \hskip0.1cm\text{is small enough})}\left(\int_M
c_n\left|\nabla
w_m\right|^2 \ dv_g+B_\epsilon\int_Mw_m^2 \ dv_g\right).
\end{eqnarray*}
Hence $$c_n
\underbrace{\left[1-\mu_2(M,g)(\mu(\mathbb{S}^n)^{-1}+\epsilon)\right]}_{>0}\int_M\left|\nabla
w_m\right|^2 \ dv_g\leq c\underbrace{\int_Mw_m^2 \ dv_g}_{\longrightarrow 0},$$
then $\int_M\left|\nabla w_m\right|^2 \ dv_g\longrightarrow 0,$ so
$\left\|w_m\right\|_{H_1^2(M)}\longrightarrow 0.$ This shows that
$w_m\longrightarrow 0$ in $H_1^2(M).$\\
We finally get that
$$1=\int_M\left|u_m\right|^{N-2}w_m^2 \ dv_g\leq \left(\int_M\left|u_m\right|^N \
dv_g\right)^\frac{N-2}{2}\underbrace{\int_Mw_m^N \ dv_g}_{\longrightarrow 0}.$$
This gives a contradiction, then $w\neq 0.$
\subsubsection{Proof of Proposition \ref{a}} Since $\la_2(g)>0$, then $$\mu_2(M,g) = \inf_{u>0} \la_2(u)\left(\int_M u^N \ dv_g\right)^\frac{2}{n} = \inf_{u\geq 0} \la_2(u)\left(\int_M u^N \ dv_g\right)^\frac{2}{n}.$$ We
mimic the proof of Lemma 3.3 in \cite{AmmannHumbert2005} by taking $w_1 = w_+ = \sup\left\lbrace
0,w\right\rbrace $ and $w_2 = w_- = \sup\left\lbrace 0, -w\right\rbrace$. This
gives that 
$$u = aw_+ + b w_-,$$
where $a,b>0$.
By Lemma 3.1, $w\in C^{2,\alpha}$, $u\in C^{0,\alpha}$ and 
Step 4 of the proof of Theorem 3.4 in \cite{AmmannHumbert2005} then shows that
$$u = \lvert w\lvert.$$
 Since $w$ is in $H_1^2(M),$ Lemma 3.1 of \cite{AmmannHumbert2005} says that $w\in L^{N + \ep}(M),$ because $w$
satisfies the equation $$L_gw =
\mu_2\lvert w\lvert^{N-2} w,$$
and standard bootstrap arguments gives that $w\in C^{3,\alpha}(M)$ for all
$\alpha<
N-2$.
\subsubsection{Proof of Proposition \ref{ip}} In this paragraph, we will see that
if $M$ is not locally conformally flat of dimension $n\geq 6,$ then
we obtain that $$\mu_2(M,g)<\mu(\mathbb{S}^n).$$ We still consider the case where
$\mu(M,g)<0$. Then there exists a positive function $v$ solution of the Yamabe
equation
\begin{equation}\label{v}
L_gv=\mu(M,g)v^{N-1}.
\end{equation}
Let $x_0$ be a point of $M$ at which the Weyl tensor is not zero (such a point exists
because the manifold is not locally conformally flat and $n\geq 4$) and $(x_1,\ldots, x_n)$ be a
system of normal coordinates at $x_0.$ For $x\in M,$ denote by $r=d(x,x_0)$ the
distance to the point $x_0.$ If $\delta$ is a small fixed number, let $\eta$ be a
cut-off function of class $C^\infty$ defined by 
$$\left\{
\begin{array}{c}
0 \leq \eta \leq 1,\\\\
\eta=1 \text{ on }B_{\delta}(x_0),\\\\
\eta=0 \text{ on } M\backslash  B_{2\delta}(x_0),\\\\
\lvert\nabla \eta \lvert\leq \frac{2}{\delta}.
\end{array}\right.$$
For all $\epsilon >0$ we define the following function
$$v_{\epsilon } =c_{\epsilon } \eta (\epsilon +r^2)^\frac{2-n}{2},$$
where $c_\epsilon$ is choosen such that $$\int_M v_{\epsilon }^N \ dv_g=1.$$
By referring to \cite{Aubin1976}
$$\lim_{\epsilon \longrightarrow 0}Y(v_{\epsilon })=\mu_1(\mathbb{S}^n),$$
where $Y(u)$ is the Yamabe functional defined by
$$Y(u)=\frac{\int_M c_n\lvert \nabla u\lvert^2+S_g u^2 \ dv_g}{\left(\int_M u^N \
dv_g\right)^\frac{2}{N}}.$$
If $(M,g)$ is not locally conformally flat, by a calculation made in
\cite{Aubin1976}, there exists a constant $C(M)>0$ such that
\begin{eqnarray}\label{Y}
   Y(v_{\ep})=\left|\;  
   \begin{matrix}
     \mu_1(\mS^n) - C(M) \ep^2 + o(\ep^2)\hfill & \hbox{ if } n > 6 \\\\
     \mu_1(\mS^n) - C(M) \ep^2 |\ln(\ep)| + o(\ep^2 |\ln(\ep)|)\hfill 
     & \hbox{ if } n = 6.
   \end{matrix}
   \right.
\end{eqnarray}
Again from \cite{Aubin1976} there exists constants $a,$ $b,$ $C_1,$ $C_2>0$, such that $$
a \ep^{\frac{n-2}{4} } \leq c_{\ep} \leq b \ep^{\frac{n-2}{4} },$$
and 
\begin{eqnarray} \label{p}
C_1 \al_{p,\ep} \leq \int_M v_{\ep}^p \,dv_g \leq C_2 \al_{p,\ep }
\end{eqnarray}
where
\[ \al_{p,\ep} = \left| \begin{array}{lll} 
\ep^{\frac{2n - (n-2) p}{4}} & \hbox{if} & p> \frac{n}{n-2};\\\\
|\ln(\ep)| \ep^{\frac{n}{4}} & \hbox{if} & p= \frac{n}{n-2};\\\\
\ep^{ \frac{(n-2) p}{4}} & \hbox{if}& p< \frac{n}{n-2}
\end{array} \right. \]
We have
\begin{eqnarray*}
\mu_2(M,g)&=&\inf_u\lambda_2(u)\left(\int_Mu^N \ dv_g\right)^{\frac{2}{n}}\\\\
&=& \inf_{\ss u \atop \ss w,w^\prime} \sup_{\lambda,\mu} \frac{\int_M L_g(\lambda w+\mu
w^\prime)(\lambda w+\mu w^\prime) \ dv_g}{\int_M u^{N-2}(\lambda w+\mu w^\prime)^2 \
dv_g}\left(\int_Mu^N \ 
dv_g\right)^{\frac{2}{n}}\\\\
&=& \inf_{\ss u \atop \ss w,w^\prime} \sup_{\lambda,\mu}F(u,\lambda w+\mu w^\prime ).
\end{eqnarray*}
Let $\lambda_\epsilon,$ $\mu_\epsilon$ such that
$$\lambda_\epsilon^2+\mu_\epsilon^2=1$$ and  $$F(v_\epsilon, \lambda_\epsilon
v+\mu_\epsilon v_\epsilon)=\sup _{(\lambda,\mu)\in \mathbb{R}^2\backslash
\left\lbrace (0,0)\right\rbrace} F(v_\epsilon, \lambda v +\mu v_\epsilon),$$
where $v$ is the function defined in the equation (\ref{v}).\\
Calculating $F(v_\epsilon, \lambda_\epsilon v+\mu_\epsilon v_\epsilon),$ we get
\begin{eqnarray*}
 F(v_\epsilon, \lambda_\epsilon v+\mu_\epsilon v_\epsilon)&=&  \frac{\int_M
L_g(\lambda_\epsilon v+\mu_\epsilon v_\epsilon )(\lambda_\epsilon  v+\mu_\epsilon
v_\epsilon ) \ dv_g}{\int_M v_\epsilon ^{N-2}(\lambda _\epsilon v+\mu_\epsilon
v_\epsilon )^2 \ dv_g}\left(\int_Mv_\epsilon ^N \ dv_g\right)^{\frac{2}{n}}\\\\
&=& \frac{\lambda_\epsilon ^2 \mu(M,g)+\mu_\epsilon ^2Y(v_\epsilon
)+2\lambda_\epsilon \mu_\epsilon \mu(M,g)\int_M v^{N-1}v_\epsilon \
dv_g}{\lambda_\epsilon ^2\int_Mv_\epsilon^{N-2}v^2 \ dv_g+
\mu_\epsilon^2+2\lambda_\epsilon \mu_\epsilon \int_M  v_\epsilon^{N-1}v \ dv_g}\\\\
&=&\frac{A_\epsilon}{B_\epsilon}.
\end{eqnarray*}
If $$\la_\ep\longrightarrow \la \neq 0, \hskip0.2cm \mu_\ep\longrightarrow \mu \neq
0,$$
then
$$F(v_\ep,\la_\ep v + \mu_\ep v_\ep)\longrightarrow \frac{\la^2 \mu(M,g) + \mu^2
\mu(\mathbb S^n)}{\mu^2}<\mu(\mathbb S^n).$$
Similarly, if $\mu = 0,$ $\la^2 = 1,$ then the numerator $A_\ep\sim \mu(M,g)< 0,$
while the denominator $B_\ep$ remains positive, which gives again that 
$$F(v_\ep,\la_\ep v + \mu_\ep v_\ep)\leq 0<\mu(\mathbb S^n),$$
which gives the desired inequality. Then, in the sequel, we assume that $\la_\ep\longrightarrow 0$ and
$\mu_\ep\longrightarrow \pm 1.$\\
 \underline{The case $n>6$}\\
Using (\ref{p}) we have
$$ \int_M v^{N-1} v_{\ep} \,dv_g \sim_{\ep \to 0} C \ep^{\frac{n-2}{4}},$$
$$\int_M v_{\ep}^{N-2} v^2 \,dv_g \sim_{\ep \to 0} C \ep,$$
and $$\int_M v_{\ep}^{N-1} v \,dv_g \sim_{\ep \to 0} C \ep^{\frac{n-2}{4}},$$
where $C$ denotes a constant that might change its value from line to line. We
distinguish two cases\\
$\bullet$ there exists a constant $C>0$ such that \begin{equation}\label{c1}
\vert \la_\ep\vert\leq C \ep^\frac{n-2}{4},\end{equation}
or\\
$\bullet$ there exists $\alpha_\ep$ such that 
\begin{equation}\label{c2}
\vert\lambda_\ep \vert= \alpha_\ep \ep^{\frac{n-2}{4}},\end{equation}
and 
$$\alpha_\ep \longrightarrow +\infty.$$ 
(possibly extracting a subsequence).
\begin{enumerate}
 \item Suppose first that (\ref{c1}) is verified. Then we have
$$\vert \la_\ep\vert\leq C \ep^\frac{n-2}{4}.$$
Hence $\la_\ep^2=O(\ep^\frac{n-2}{2}),$ so
$\mu^2_\ep=1-\la^2_\ep=1+O(\ep^\frac{n-2}{2}).$
Therefore $$\mu_\ep=1+O(\ep^\frac{n-2}{2}).$$
This gives
\begin{eqnarray*}
A_\epsilon &=& O(\ep^{\frac{n-2}{2}}) + ( 1 +O(\ep^{\frac{n-2}{2}})) \Big( \mu
(\mathbb{S}^n) -C(M)\ep^2 +o(\ep^2)\Big) + O(\ep^{\frac{n-2}{2}})\\
&=& \mu (\mathbb{S}^n)-C(M)\ep^2 + O(\ep^{\frac{n-2}{2}}) +o(\ep^2).
\end{eqnarray*}
Since $\frac{n-2}{2} > 2$,
$$A_\ep =\mu (\mathbb{S}^n) -C(M)\ep^2 +o(\ep^2),$$
and 
$$B_\ep = O(\ep^{\frac{n-2}{2}+1}) + 1 + O(\ep^{\frac{n-2}{2}})+
O(\ep^{\frac{n-2}{2}}) = 1+ o(\ep^2).$$
Then, 
$$\frac{A_\ep}{B_\ep} = \mu(\mathbb{S}^n)- C(M) \ep^2 + o(\ep^2) < \mu(\mathbb{S}^n).$$
\item Assume now that (\ref{c2}) is fulfilled. In this case
\begin{eqnarray*}
\frac{A_\ep}{B_\ep} &=& \frac{\lambda^2_\ep \mu (M, g) + (1 - \lambda^2_\ep)
Y(v_\ep) + \lambda_\ep O(\ep^{\frac{n-2}{4}})} { \lambda^2_\ep
O(\ep) + (1- \lambda^2_\ep) + 2\lambda_\ep \mu_\ep O(\ep^{\frac{n-2}{4}})}\\
&=& \frac{\lambda^2_\ep \mu (M, g) + (1 - \lambda^2_\ep) Y(v_\ep) + 
o(\la^2_\ep)}{o(\la^2_\ep) + (1-\la^2_\ep) + o(\la^2_\ep)}\\
&=& \frac{\la^2_\ep \mu(M, g)}{1-\la_\ep^2 + o(\la^2_\ep)} + \frac{Y(v_\ep)}{ 1 +
\frac{o(\la_\ep^2)}{\mu^2_\ep}} + o(\la_\ep^2)\\\\
&\leq& \mu(M,g) \la_\ep^2 + \mu(\mathbb{S}^n)(1 + o(\la_\ep^2)) + o(\la_\ep^2)\\
&\leq& \mu(\mathbb{S}^n) + \mu(M,g) \la_\ep^2 + o(\la_\ep^2)\\
&<& \mu(\mathbb{S}^n),
\end{eqnarray*}
because $\mu(M,g)<0$ and $Y(v_\ep) \leq \mu(\mathbb{S}^n).$
\end{enumerate}
\underline{The case $n =6$}\\
Since \begin{eqnarray*}
 \int_M v_\ep^{N-2} v^2 dv_g &\sim_{\ep \rightarrow 0}& C \ep,\\
\int_M v^{N-1}v_\ep dv_g &\sim_{\ep \rightarrow 0}& C \ep,\\
\int_M v_\ep^{N-1} v dv_g &\sim_{\ep \rightarrow 0}& C \ep,
\end{eqnarray*}
then
$$A_\ep = \la_\ep^2 \mu(M,g) + \mu^2_\ep Y(v_\ep) + 2\la_\ep \mu_\ep O(\ep),$$
$$B_\ep = \la^2_\ep O(\ep) + \mu^2_\ep + 2\la_\ep \mu_\ep O(\ep).$$
Again, we have two cases to study
\begin{enumerate}
\item If $\lvert \la_\ep\lvert\leq C\ep,$ then $$\la_\ep^2\leq C\ep^2.$$  This implies
$$A_\ep = \mu(\mathbb{S}^n) - C\ep^2\lvert \ln(\ep)\lvert + o(\ep^2\lvert
\ln(\ep)\lvert)$$
and 
$$B_\ep = 1+O(\ep^2) = 1 + o(\ep^2 \lvert \ln(\ep)\lvert).$$
Hence
$$\frac{A_\ep}{B_\ep}< \mu(\mathbb{S}^n).$$
\item If $\lvert \la_\ep\lvert = \alpha_\ep \ep,$ with $\alpha_\ep \longrightarrow
+\infty.$ Since $Y(v_\ep) \leq \mu(\mathbb{S}^n),$ therefore
$$A_\ep = \alpha_\ep^2 \ep^2 \mu(M,g) + \mu^2_\ep \mu(\mathbb{S}^n) + o(\alpha_\ep^2
\ep^2),$$
and $$B_\ep = \mu^2_\ep + o(\alpha_\ep^2 \ep^2).$$
Therefore
\begin{eqnarray*}\frac{A_\ep}{B_\ep} &=& \mu(M,g) \frac{\alpha_\ep^2 \ep^2}{1+o(1)}
+ \frac{\mu(\mathbb{S}^n)}{1 + o(\alpha_\ep^2 \ep^2)} + \frac{o(\alpha_\ep^2 \ep^2)}{1
+ o(1)}\\
 &=& \mu(M,g)\alpha_\ep^2 \ep^2 + \mu(\mathbb{S}^n) + o(\alpha_\ep^2 \ep^2)\\
&<& \mu(\mathbb{S}^n).
\end{eqnarray*}
\end{enumerate}
This ends the proof of Propositon \ref{ip}.
\hfill $\square$\\
So we get a solution $w$ having a changing sign of the equation 
$$L_g w = \mu_2 \lvert w\lvert ^{N-2}w.$$
Finally, to obtain the resultat announced in Theorem \ref{th1}, it suffices to
set $$w^\prime = \mu_2^\frac{n-2}{4} w,$$
then $w^\prime$ verifies 
$$L_gw^\prime = \ep\lvert w^\prime\lvert^{N-2} w^\prime,$$
with $\ep = 1 = \hbox{ sign}(\la_2(g))$.
\section{The case $\la_2<0$}
In this section, we will show that in all cases, there exists a nodal solution of
the equation
$$L_g w = C_0 \lvert w\lvert^{N-2}w,$$
where $C_0$ is a negative constant.\\

\noindent First, since $\mu<0,$ we assume in the whole section
that the metric $g$ is such that $S_g=-1.$
In this context, the approach will be different. Indeed, the second Yamabe invariant 
is not well defined as shown in the following proposition:
\begin{prop}\label{-infty}
Let $M$ be a compact Riemannian manifold of dimension $n\geq 3$. Suppose
that $\la_2< 0$, then
$$\inf_u \la_2(u)\left(\int_M u^N \ dv_g\right)^\frac{2}{n} = -\infty.$$
\end{prop}
The proof will be detailed in Subsection \ref{inf=-infty}.\\
We will use a new functional $$I_g(u)= \frac{\left(\int_M\lvert L_g
u\lvert^\frac{2n}{n+2} \ dv_g\right)^\frac{n+2}{n}}{\lvert\int_M uL_gu \ dv_g\lvert}.$$
We study $\alpha:= \inf I_g(u)$ where the infimum is taken over the functions 
$u\in H_2^\frac{2n}{n+2}(M)$ such that
$$\int_M uL_gu \ dv_g <0,$$
and with the following constraint
$$\int_M \lvert u \lvert^{N-2} u \ v \ dv_g=0,$$
for any function $v\in \ker L_g.$\\

\noindent We will show that $\alpha$ is a conformal invariant. We obtain also that the infimum of this
functional is attained by a function $u.$ We set 
$$v = \lvert L_gu\lvert^\frac{-4}{n+2} L_gu,$$
and we will observe that $v$  has the following properties:\\
$\bullet$ $v$ is a solution of the equation $$L_g v = \alpha^\prime \lvert
v\lvert^{N-2}v,$$ where $\alpha^\prime<0$ (i.e. has same sign than
$\la_2$).\\
$\bullet$ $v$ has a changing sign.\\
$\bullet$ $v$ is of class $C^{3,\alpha}(M)$ ($\alpha < N-2$).
\subsubsection{Conformal invariance of $\alpha$}
Let $\widetilde g = \phi^\frac{4}{n-2} g$ be a conformal metric, $\phi$
a smooth positive function. Then
$$dv_{\widetilde g} = \phi^\frac{2n}{n-2} dv_g,$$
and
$$L_{\widetilde g}u = \phi^{-\frac{n+2}{n-2}}L_g(u\phi),$$
for all functions $u.$
\begin{enumerate}
 \item Remark that $I_{\widetilde g}(u) = I_g(u\phi).$
\begin{eqnarray*}
 I_{\widetilde g} (u)&=& \frac{\left(\int_M \lvert L_{\widetilde g}
u\lvert^\frac{2n}{n+2} \ dv_{\widetilde g}\right)^\frac{n+2}{n}}{\lvert \int_M u
L_{\widetilde g}u \ dv_{\widetilde g}\lvert}\\\\
&=& \frac{\left(\int_M \lvert \phi\lvert^\frac{-2n}{n-2} \lvert
L_g(u\phi)\lvert^\frac{2n}{n+2} \phi^\frac{2n}{n-2} \
dv_g\right)^\frac{n+2}{n}}{\lvert \int_M u\phi^\frac{-(n+2)}{n-2} L_g (u\phi)
\phi^\frac{2n}{n-2} \ dv_g\lvert}\\\\
&=& \frac{\left(\int_M \lvert L_g(u\phi)\lvert ^\frac{2n}{n-2} \
dv_g\right)^\frac{n+2}{n}}{\lvert \int_M u\phi L_g(u\phi) \ dv_g\lvert}\\\\
&=& I_g(u\phi),
\end{eqnarray*}
where we have used
\begin{eqnarray*}
 \int_M uL_{\widetilde g}u \ dv_{\widetilde g} &=& \int_M u\phi^{-\frac{n+2}{n-2}}
L_g(u\phi)\phi^\frac{2n}{n-2} \ dv_g\\
&=& \int_M (u\phi)L_g(u\phi) \ dv_g.
\end{eqnarray*}
\item Assume that for any $v\in \ker L_{\widetilde g},$ we have
$$\int_M \lvert u\lvert^{N-2} uv \ dv_{\widetilde g} = 0.$$
Then, for any $v^\prime \in \ker L_g,$ we obtain
$$\int_M \lvert u\phi\lvert^{N-2} (u\phi)v^\prime \ dv_g = \int_M \lvert
u\lvert^{N-2} u (v^\prime \phi^{-1}) \ dv_{\widetilde g} = 0,$$
since $$L_{\widetilde g}(v^\prime \phi^{-1}) = \phi^{-\frac{n+2}{n-2}}
L_g(v^\prime) = 0,$$
i.e. $$v^\prime \phi^{-1}\in Ker L_{\widetilde g}.$$ \hfill$\square$
\end{enumerate}
\subsubsection{Proof of Proposition \ref{-infty}}\label{inf=-infty}
Assume that $\la_2(g)< 0,$ and choose $u>0.$\\
By Lemma \ref{existancedevetw}, there exists two functions $v_1$ and $v_2$
solutions of the following equations
$$L_g v_1 = \la_1(u) \lvert u\lvert^{N-2} v_1,$$ 
and
$$L_g v_2 = \la_2(u) \lvert u\lvert^{N-2} v_2,$$ 
such that
$$\int_M \lvert u\lvert^{N-2} v_1v_2 \ dv_g=0.$$
Let $v_\ep$ the function defined in Section \ref{section}, and let
$V=\left\lbrace v_1,v_2\right\rbrace.$ For all $v\in V,$ we get 
$$\lim_{\ep\longrightarrow 0}\int_Mv_\ep^{N-2} v^2 \ dv_g = 0.$$ 
Since $\la_1(u)<0$ and $\la_2(u)<0,$ then for $\ep$ sufficiently small, we have
$$\lim_{\ep\to 0}\left(\sup_{v\in V}\frac{\int_M (L_gv)(v) \ dv_g}{\int_M v_\ep^{N-2}v^2 \ dv_g}\right) = 
-\infty,$$
hence $$\lim_{\ep\to 0}\left(\inf_u \la_2(u)\left(\int_M u^N \ dv_g\right)^\frac{2}{n}\right) = 
 -\infty.$$ \hfill$\square$
\subsubsection{The infimum of the functional $I_g$ is attained} Let $(u_m)_m$ be a
minimizing sequence, i.e., $$\lim_{m\longrightarrow \infty}
I_g(u_m) = \alpha,$$
 with 
$$\int_M \lvert u_m \lvert^{N-2} \ u_m \ v \ dv_g=0, \hskip0.2cm\forall \ v\in \ker
L_g.$$
We can assume that 
\begin{equation}\label{-1}
\int_M u_mL_gu_m \ dv_g = -1.
\end{equation}
Then
$$
\alpha = \lim_{m\longrightarrow\infty} \left(\int_M\lvert L_gu_m\lvert^\frac{2n}{n+2} \
dv_g\right)^\frac{n+2}{n}.$$
Now we show that $(u_m)_m$ is a bounded sequence in $H_2^\frac{2n}{n+2}(M).$\\
We proceed by contradiction and we assume that, up to a subsequence, $\lim \|
u_m\|_{H_2^\frac{2n}{n+2}(M)}= +\infty$. Let $$v_m = \frac{u_m}{\|
u_m\|_{H_2^\frac{2n}{n+2}(M)}}.$$
Since $\| v_m\|_{H_2^\frac{2n}{n+2}(M)} = 1,$ $(v_m)_m$ is a bounded sequence
in $H_2^\frac{2n}{n+2}(M),$ and therefore there exists $v \in H_2^\frac{2n}{n+2}(M)$
such that after restriction to a subsequence
$$v_m\rightharpoonup v\text{ in }H_2^\frac{2n}{n+2}(M),$$
$$v_m\longrightarrow v\text{ in }L^2(M).$$
By standard arguments, we get
$$\left(\int_M\lvert L_gv\lvert^\frac{2n}{n+2} \ dv_g\right)^\frac{n+2}{n}\leq
\underbrace{\liminf_m \left(\int_M\lvert L_gv_m\lvert^\frac{2n}{n+2} \
dv_g\right)^\frac{n+2}{n}}_{=0}.$$
This gives 
$$L_gv = 0,$$
hence $$v\in \ker L_g.$$
We have for all function $v^\prime \in \ker L_g,$
$$\int_M \lvert v_m\lvert^{N-2}v_m v^\prime \ dv_g = \frac{\int_M \lvert
u_m\lvert^{N-2}u_m v^\prime \ dv_g}{\| u_m\|^{N-1}_{H_2^\frac{2n}{n+2}}} = 0.$$
In particular for $v^\prime = v$, 
$$\int_M \lvert v_m\lvert^{N-2}v_m v \ dv_g = 0\longrightarrow_{m\to \infty}\int_M
v^N \ dv_g,$$
so
\begin{equation}\label{0}
v=0.\end{equation}
According to the regularity Theorem 3.75 in \cite{Aubin98}, we have
$$1 = \| v_m\|_{H_2^\frac{2n}{n+2}}\leq C\left[
\underbrace{\|L_gv_m\|_{L^\frac{2n}{n+2}}}_{\longrightarrow 0} +
\|v_m\|_{L^\frac{2n}{n+2}}\right].$$
Passing to the limit, we obtain 
$$\int_M v^\frac{2n}{n+2} \ dv_g\geq \frac{1}{C},$$
which gives a contradiction. We deduce that $(u_m)_m$
is a bounded sequence in $H_2^\frac{2n}{n+2}(M).$ Then, after restriction to a
subsequence, there exists $u$ in $H_2^\frac{2n}{n+2}(M)$ such that
$$u_m\rightharpoonup u \text{ in } H_2^\frac{2n}{n+2}(M),$$
$$u_m\rightharpoonup u  \text{ in }H_1^2(M),$$
$$u_m\longrightarrow u \text{ in } L^2(M).$$
Further, we have $$\left(\int_M \lvert L_gu\lvert^\frac{2n}{n+2} \
dv_g\right)^\frac{n+2}{n}\leq \liminf_m \left(\int_M \lvert
L_gu_m\lvert^\frac{2n}{n+2} \
dv_g\right)^\frac{n+2}{n}.$$
Moreover,
\begin{align}
 \int_M uL_gu \ dv_g =& \int_M \lvert \nabla u\lvert ^2 \ dv_g - \int
_M u^2 \
dv_g \nonumber\\
\leq& \liminf _m  \int_M \lvert \nabla u_m\lvert ^2 \ dv_g - \int _M u_m^2 \ dv_g
\nonumber\\
=& \liminf _m  \int_M u_mL_gu_m \ dv_g = -1.\label{liminf}
\end{align}
Therefore $$ \int_M uL_gu \ dv_g < 0,$$
and $$u\neq 0.$$
Finally, with (\ref{liminf})
\begin{eqnarray*}
 I_g(u) &=& \frac {\left(\int_M \lvert L_gu\lvert^\frac{2n}{n+2} \ dv_g
\right)^\frac{n+2}{n}}{\lvert \int_M uL_gu \ dv_g\lvert}\\
&\leq& \liminf_m \frac {\left(\int_M \lvert L_gu_m\lvert^\frac{2n}{n+2} \ dv_g
\right)^\frac{n+2}{n}}{\lvert \int_M u_mL_gu_m \ dv_g\lvert}= \liminf_m I_g(u_m) =
\alpha.
\end{eqnarray*}
Hence the result is proved i.e. $I_g(u) = \alpha$.\\
\underline{Euler equation}\\
Notice that
 $$\int_M \lvert u\lvert^{N-2} \ u \ v^\prime \ dv_g = 0, \text{ for any
function }v^\prime \in \ker L_g.$$ 
In particular, $\al \not =0$. Remark also that 
$$\int_M uL_gu \ dv_g = -1.$$
Indeed, the relation $\int_M uL_gu \ dv_g < -1$ would imply that $I_g(u)<  \lim I_g(u_m) = \al$. 
We now write Euler equation of $u$.
Let $\left\lbrace u_1, \ldots ,u_k\right\rbrace$ be a base of $\ker L_g.$ By the
Lagrange multipliers theorem, there exists real numbers $ \la_1,\ldots ,\la_k$
for which, for all function $\phi\in C^\infty(M),$ we get
$$\frac{d}{dt}|_{t=0}I_g(u+t\phi) = \Sigma_i \la_i\frac{d}{dt}|_{t=0} g_i(u+t\phi),$$
where $$g_i(u) = \int_M \lvert u \lvert^{N-2} \ u \ u_i \ dv_g.$$
Setting $a = \left(\int_M \lvert L_gu\lvert^\frac{2n}{n+2} \
dv_g\right)^\frac{n+2}{n}$, one checks $$2 a^\frac{2}{n+2}\int_M
\lvert L_gu\lvert^\frac{-4}{n+2}L_gu L_g\phi\ dv_g + 2 a\int_M \phi L_gu \
dv_g = (N-1) \Sigma_i
\la_i\int_M \lvert u \lvert^{N-2} \ \phi \ u_i \ dv_g.$$
If $\phi\in \ker L_g,$ this last equation implies that
$$\Sigma_i \la_i\int_M \lvert u \lvert^{N-2} \ \phi \ u_i \ dv_g = 0.$$
Then, for $\phi = \Sigma_i\la_i \ u_i \in \ker L_g,$ we have 
$$\int_M \lvert u\lvert^{N-2} \ \phi^2 \ dv_g = 0.$$
Therefore 
\begin{eqnarray*}
 \lvert u\lvert^{N-2} \ \phi^2 = 0&\Rightarrow & \lvert u\lvert^{N-2} \ \phi = 0\\
&\Rightarrow &\Sigma_i \ \la_i \ \lvert u\lvert^{N-2} \ u_i = 0.
\end{eqnarray*}
This gives, for any function $\phi$ (in $\ker L_g$ or not), that
$$\Sigma_i \la_i\int_M \lvert u\lvert^{N-2} \phi \ u_i \ dv_g = 0.$$
Then $u$ verifies in the sense of distributions the following equation 
\begin{equation} \label{u}
L_g\left(\lvert L_gu\lvert ^\frac{-4}{n+2} L_gu\right) = \alpha^\prime L_gu,
\end{equation}
where $$\alpha^\prime = - \alpha ^\frac{n}{n+2}= -\int_M \lvert
L_gu\lvert^\frac{2n}{n+2} \ dv_g.$$
We set
$$v = \lvert L_gu\lvert ^\frac{-4}{n+2} L_gu \in L^N(M),$$
then $$\lvert v\lvert = \lvert L_gu\lvert^{1-\frac{4}{n+2}} = \lvert
L_gu\lvert^\frac{n-2}{n+2}.$$
Hence,
$$ L_gu = \lvert v\lvert ^{N-2} \ v.$$
Replacing each term by its value in Equation (\ref{u}), we obtain
$$L_gv = \alpha ^\prime \lvert v\lvert ^{N-2} \ v.$$
\underline{Regularity of $v$}\\
We have $u\in H_2^\frac{2n}{n+2}(M),$ then $L_gu \in L^\frac{2n}{n+2}(M).$
Therefore $$v\in L^N(M),$$
since $\lvert v\lvert ^N = \lvert L_g u\lvert^\frac{2n}{n+2}.$
Moreover, in the sense of distributions 
\begin{equation}\label{prime}
L_gv = \alpha ^\prime \lvert v\lvert ^{N-2} \ v,
\end{equation}
this implies that 
$$\lvert L_g v \lvert = \lvert \alpha ^\prime \lvert \lvert
v\lvert^{N-1},$$
hence $L_gv \in L^\frac{N}{N-1} (M)= L^\frac{2n}{n+2}(M),$ therefore $v \in
H_2^\frac{2n}{n+2}(M)\subset H_1^2(M).$\\
Using Lemma 3.1 of \cite {AmmannHumbert2005}, we get 
$$v\in L^{N+\ep}(M),$$
By a standard bootstrap argument, we show that $v\in C^{3,\alpha}(M) (\alpha <
N-2).$\\
Calculating now $I_g(v),$ using (\ref{prime}), we have
\begin{eqnarray*}
 I_g(v) &=& \frac{\left(\int_M\lvert L_g v\lvert^\frac{2n}{n+2} \
dv_g\right)^\frac{n+2}{n}}{\lvert\int_M vL_gv \ dv_g\lvert}\\\\
&=& \frac{\alpha^{\prime 2}\left(\int_M \lvert v\lvert^{(N-1)\times \frac{2n}{n+2}} \
dv_g\right)^\frac{n+2}{n}}{\lvert \alpha^\prime\lvert \int_M \lvert v\lvert^N \
dv_g}\\\\
&=& \alpha^\frac{n}{n+2} \frac{\left(\int_M\lvert v \lvert^\frac{2n}{n-2} \
dv_g\right)^\frac{n+2}{n}}{\int_M \lvert v\lvert^\frac{2n}{n-2} \ dv_g}\\\\
&=& \alpha^\frac{n}{n+2} \left(\int_M \lvert v\lvert^\frac{2n}{n-2} \
dv_g\right)^\frac{2}{n}\\\\
&=& \alpha^\frac{n}{n+2} \left(\int_M \lvert L_gu\lvert^\frac{2n}{n+2} \
dv_g\right)^\frac{2}{n}\\\\
&=& \alpha^\frac{n}{n+2} \alpha^\frac{2}{n+2} = \alpha.
\end{eqnarray*}
The function $v$ satisfies that, for any function $v^\prime \in\ker L_g,$
$$\int_M \lvert v\lvert ^{N-2}v v^\prime \ dv_g = 0.$$ 
Indeed,
\begin{eqnarray*}
 \int_M \lvert v\lvert^{N-2} \ v \ v^\prime \ dv_g &=& \int _M L_gu \ v^\prime \ dv_g\\
&=& \int_M u \ L_g v^\prime \ dv_g \\
&=&0.
\end{eqnarray*}
\underline{$\bullet$ $v$ has changing sign}\\
We proceed by contradiction and assume that $v\geq 0$. Since $v \not=0$, we deduce from the maximum principle that $v>0$. In addition, Equation (\ref{prime}) says that there exists an $i$ such that $\alpha^\prime
= \la_i(v)$. The only positive eigenfunctions are the ones associated to $\lambda_1$ and hence $\al'=\lambda_1(v)$. By Proposition \ref{existancedevetw}, there exists a function $w$ solution
of the following equation 
$$L_gw = \la_2(v) \lvert v\lvert^{N-2} w.$$
\begin{eqnarray*}
 I_g(w) &=& \frac{\left(\int_M \lvert L_gw\lvert^\frac{2n}{n+2} \
dv_g\right)^\frac{n+2}{n}}{\lvert \int_M wL_gw \ dv_g\lvert}\\\\
&=& \frac{\lvert \la_2(v)  \lvert ^2 \left(\int_M \lvert v\lvert^{(N-2)\times
\frac{2n}{n+2}} \lvert w\lvert^\frac{2n}{n+2} \ dv_g\right)^\frac{n+2}{n}}{\lvert
\la_2(v)
 \lvert \int_M \lvert v\lvert^{N-2} \ w^2 \ dv_g}.
\end{eqnarray*}
By applying the H\"older inequality with $p = \frac{n+2}{n}$ and $q = \frac{n+2}{2},$ we get
\begin{eqnarray*}
 \int_M \lvert v\lvert^{(N-2)\times \frac{2n}{n+2}} \ \lvert w\lvert^\frac{2n}{n+2}
\ dv_g &=& \int_M \lvert v\lvert^{(N-2)\times \frac{n}{n+2}} \ \lvert
w\lvert^\frac{2n}{n+2} \  \lvert v\lvert^{(N-2)\times \frac{n}{n+2}} \ dv_g\\\\
&\leq& \left(\int_M \lvert v\lvert^{N-2} \ w^2 \ dv_g\right)^\frac{n}{n+2}
\left(\int_M \lvert
v\lvert^{\frac{4}{n-2}\times\frac{n}{2}} \ dv_g\right)^\frac{2}{n+2}.
\end{eqnarray*}
Therefore 
\begin{eqnarray*}
I_g(w) &\leq& \lvert \la_2(v)  \lvert \left(\int_M \lvert v \lvert^\frac{2n}{n-2} \
dv_g\right)^\frac{2}{n}\\
&=& \lvert \la_1(v)  \lvert \left(\int_M \lvert L_gu \lvert^\frac{2n}{n+2} \
dv_g\right)^\frac{2}{n}\\
&=& \alpha^\frac{n}{n+2}\alpha^\frac{2}{n+2} = \alpha,
\end{eqnarray*}
since by assumption $\la_1(v) = \alpha^\prime = \alpha^\frac{n}{n+2},$ which gives
a contradiction.\\
Then $v$ is a nodal solution of the equation
$$L_g v = \alpha^{\prime} \lvert v\lvert^{N-2} v,$$
where $\alpha^{\prime} <0$.
Setting 
$$v^\prime: = \lvert \alpha\lvert^\frac{n-2}{4},$$
we obtain that $v^\prime$ is a solution
of the  equation
$$L_g v^\prime =  \ep \lvert v^\prime\lvert^{N-2} v^\prime$$
with $\ep = -1 \hbox{ sign }(\lambda_2(g))$. This ends the proof of Theorem \ref{th1}.
\bibliographystyle{amsalpha}
\bibliography{biblio}

\end{document}